\documentclass[12pt,twoside]{article}
\usepackage{graphicx}
\newtheorem{theorem}{\bf Theorem}[section]

\newtheorem{lemma}[theorem]{\bf Lemma}
\newtheorem{corollary}[theorem]{\bf Corollary}

\input{amssym}

\date{September 10, 2009}
\begin{document}
\title{{\Large Large stars with few colors}}
\author{\small A. Khamseh$^{\textrm{a}}$, G.R. Omidi$^{\textrm{a},\textrm{b},1}$\\
\small  $^{\textrm{a}}$Department of Mathematical Sciences, Isfahan
University
of Technology,\\ \small Isfahan, 84156-83111, Iran\\
\small  $^{\textrm{b}}$School of Mathematics, Institute for Research
in Fundamental Sciences (IPM),\\
\small  P.O.Box:19395-5746, Tehran,
Iran\\
\small \texttt{E-mails:khamseh@math.iut.ac.ir, romidi@cc.iut.ac.ir}}
\date {}

\maketitle\footnotetext[1] {\tt This research was in part supported
by a grant from IPM (No.88050012)} \vspace*{-0.5cm}

\begin{abstract}

A  recent question in generalized Ramsey theory is that
for fixed positive integers $s\leq t$,
at least how many vertices can be covered by the vertices
of no more than $s$ monochromatic members of the family
$\cal F$ in every edge coloring of
$K_n$ with $t$ colors. This is related to
an old problem of
Chung and Liu: for graph $G$ and integers
$1\leq s<t$ what is the smallest positive integer $n=R_{s,t}(G)$ such
that every coloring of the edges of $K_n$ with $t$ colors contains a copy of $G$
with at most $s$ colors.
We answer this question
when $G$ is a star and $s$ is either $t-1$ or $t-2$  generalizing
the well-known result of  Burr and Roberts.

\end{abstract}
\section{Introduction}

Ramsey theory is an area
of combinatorics  which uses techniques from many branches of
mathematics and is currently among the most active areas in
combinatorics. Let $G_1, \ldots , G_c$ be graphs. The{\it{ Ramsey
number}} denoted by $r(G_1, \ldots , G_c)$ is defined to be the
least number $p$ such that if the edges of the complete graph
$K_p$ are arbitrarily colored with $c$ colors, then for some $i$
the spanning subgraph whose edges are colored with the $i$-th
color contains  $G_i$. More information about the Ramsey numbers
of known graphs can be found in the survey \cite{sur}.

\bigskip
There are
various types of Ramsey numbers that are important in the study of
classical Ramsey numbers and also hypergraph Ramsey numbers.
A question  recently proposed
by  Gy\'{a}rf\'{a}s et al. in \cite{gyarfas};
for fixed positive integers $s\leq t$,
at least how many vertices can be covered by the vertices
of no more than $s$ monochromatic members of the family
$\cal F$ in every edge coloring of
$K_n$ with $t$ colors.
This is  related to an old problem
of Chung and Liu \cite{chli2}:
for a given graph $G$ and for fixed  $1\leq s<t$, find the smallest $n=R_{s,t}(G)$ such that in
every $t$-coloring of the edges of $K_n$ there is a copy of $G$ colored with at most $s$
colors. Note that for $s=1$ this is the same Ramsey number $r_t(G)$.
Several problems and interesting
conjectures was presented in \cite{gyarfas}.
A basic problem here is to find
the largest $s$-colored element of $\cal F$
that can be found in every $t$-coloring of $K_n$.
The answer for  matchings when $s=t-1$ was given in \cite{gyarfas};
every $t$-coloring of $K_n$ contains a $(t-1)$-colored matching
of size $k$ provided that $n\geq 2k+[\frac{k-1}{2^{t-1}-1}]$.
Note that for $t=2,3,4$, we can guarantee the existence
of a $(t-1)$-colored path on ${2k}$ vertices instead of  a matching
of size $k$.
This was proved  in \cite{gerencser}, \cite{meenakshi}, and
\cite{khomidiforests}, respectively.
For complete graphs the problem was partially
answered  in \cite{chli2} and \cite{harborth}.
Naturally,  for these graphs the answer is very few known
and there are many open problems.
For stars, when  $s=1$ it is the  well-know
result  of  Burr and Roberts \cite{ramseyforstars}, and when $s=t-1=2$
it was determined in   \cite{chli1}.

\bigskip
In this paper we find the value of $R_{s,t}(G)$
when $G$ is a star and $s$ is either $t-1$ or $t-2$. This will generalize
the  results of \cite{ramseyforstars} and   \cite{chli1}.
The paper is organized as follows.
In section \ref{more}, we  give the upper bound and
lower bound of $R_{t-l, t}(K_{1, n})$ for given integer $l\geq 1$. In sections
\ref{secdelta} and \ref{l-2}, we determine the values of $R_{t-1,
t}(K_{1, n})$ and  $R_{t-2, t}(K_{1, n})$, respectively. As usual,
we only concerned with undirected simple finite
graphs and for the vertex $v$ of $G$ the set of edges adjacent to
$v$ in $G$ is denoted by  $E_G(v)$.


\section{\normalsize{Some bounds}}\label{more}

In this section, we find some bounds for $R_{t-l, t}(K_{1, n})$.
The {\it{Tur$\acute{a}$n number}} $ex(H, p)$ is the maximum number
of edges in a graph on $p$ vertices which is $H$-free, i.e. it
does not have $H$ as a subgraph. It is easily seen that $ex(K_{1,
n}, p)\leq \frac{p(n-1)}{2}$. This fact yields an upper bound for
$R_{t-l, t}(K_{1, n})$ as we see in the following theorem.

\begin{theorem}\label{thm1}
Suppose that $t'=[t/l]$, then $R_{t-l, t}(K_{1, n})\leq p$ for
$p>\frac{t'n-1}{t'-1}$.
\end{theorem}
{\it Proof.} Consider an edge coloring of  $K_p$ with $t$ colors.
Divide these $t$ colors into $t'=[t/l]$ classes each of which
contains $l$ colors except the last one which may contains more
colors. There exist $l$ colors  with at most $\left[\frac{1}{t'}
{p \choose 2}\right]$ edges. Thus the remaining $t-l$ colors
appear on at least ${p \choose 2}-\left[\frac{1}{t'} {p \choose
2}\right]$ edges and the existence of $K_{1, n}$ with these $t-l$
colors is guaranteed if
$${p \choose 2}-\left[\frac{1}{t'} {p \choose
2}\right]>\frac{p(n-1)}{2}.$$ So if $p>\frac{t'n-1}{t'-1}$, the
above inequality is fulfilled and there
exists a $K_{1,n}$ with at most $t-l$ colors.$\hfill \dashv $ \\

The next theorem gives a lower bound for $R_{t-l, t}(K_{1, n})$.

\begin{theorem}\label{thm2}
Let $y=\left[\frac{t(n-l+1)-l}{t-l}\right]$. Then $R_{t-l,
t}(K_{1, n})>y-\epsilon$ where $\epsilon=1$ if $y$ is odd and
$\epsilon=0$, otherwise.

\end{theorem}
{\it Proof.} Let $p=y-\epsilon$. It is sufficient to give an edge
coloring of $K_p$ such that the set of colors appear on the edges of
every $K_{1, n}$ contains at least $t-l+1$ colors. By Vizing's
theorem, there exists a proper edge coloring of $K_p$ with $p-1$
colors. Let $p-1=qt+r$, $0\leq r\leq t-1$. We partition the above
$p-1$ colors into $t$ classes each of which contains
$q=\left[\frac{p-1}{t}\right]$ colors except the last one which may
contains $(p-1)-q(t-1)$ colors. Every $K_{1, n}$ contains at least
$t-l+1$ colors if
$$n>(t-l-1)q+p-1-(t-1)q=(p-1)-lq.$$
The above inequality holds if $\frac{p-1}{t} \geq \frac{p-n-1}{l}+1$
or equivalently, $p\leq \frac{t(n-l+1)-l}{t-l}$ as asserted in
Theorem \ref{thm2}. So there is no $K_{1, n}$
with at most $t-l$ colors, that is, $R_{t-l, t}(K_{1, n})> p$.$\hfill \dashv $ \\

Combining Theorems \ref{thm1} and \ref{thm2}, we have an
approximation of the value of $R_{t-l, t}(K_{1, n})$. For the small
values of $l$ this approximation is closer to the exact value. In
particular, for $l=1, 2$, we have the following corollaries.

\begin{corollary}\label{cor1}
Let $x=\left[\frac{nt-1}{t-1}\right]$. Then
$$x\leq R_{t-1, t}(K_{1, n}) \leq x+1.$$ In particular, when
$x$ is even, then $R_{t-1, t}(K_{1, n})=x+1$.
\end{corollary}

\begin{corollary}\label{thm3}
Let $t\geq 4$, $t'=[t/2]$ and $x=\left[\frac{nt'-1}{t'-1}\right]$.
Then
$$x-2\leq R_{t-2, t}(K_{1, n}) \leq x+1.$$ In particular, when
$\left[\frac{t(n-1)-2}{t-2}\right]$ is even, then $x-1\leq R_{t-2,
t}(K_{1, n}) \leq x+1$.
\end{corollary}

\vspace{.5cm} \noindent {\bf Remark}. Let  $v_1, \ldots, v_x$ be
vertices of $K_x$, where $x$ is odd. Eliminating $v_x$, there
exists corresponding matching $M_{v_x}$ 
containing $(x-1)/2$ parallel edges $v_1v_{x-1}, v_2v_{x-2}
,\ldots, v_{(x-1)/2}v_{(x+1)/2}$. Order these edges as above.
Similarly, for each vertex $v_i$, $1\leq i\leq x-1$, there exists
the matching $M_{v_i}$ containing $(x-1)/2$ ordered edges. These
matchings are used to construct certain edge colorings of $K_x$,
for example in the proof of following key lemmas.

\begin{figure}[h]
  \begin{center}
  \includegraphics[width=9cm]{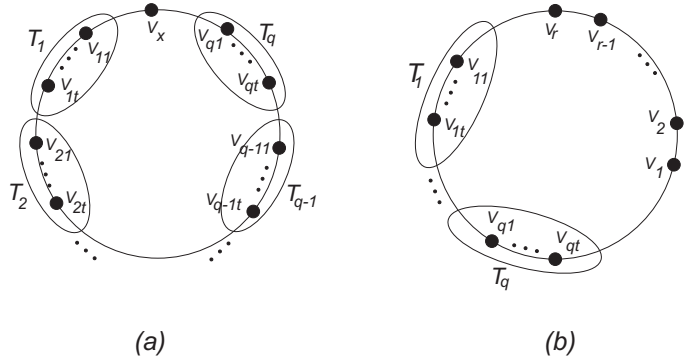}
  \caption{Partitions of the vertices of $K_x$ }\label{G1}
  \end{center}
\end{figure}

\begin{lemma}\label{lemmaA}
Suppose that $q$ is even and $x-1=tq$. There exists an edge
coloring of $K_x$ with $t$ colors such that  the set of all
neighbors of every vertex contains $q$ edges of any  color.
\end{lemma}
{\it Proof.} Partition the  vertices of $K_x$ as a single vertex
$v_x$ plus $q$ classes $T_1, \ldots, T_q$ where  $T_i$  contains $t$
vertices say $v_{i1}, \ldots, v_{it}$. Set $q/2$ classes $T_1,
\ldots, T_{q/2}$ on one side of $v_x$ and $q/2$ classes $T_{q/2+1},
\ldots, T_q$ on the other side of $v_x$ (see $(a)$ of figure
\ref{G1}). For each vertex $v_{ij}$, $1\leq j \leq t$ and $1\leq
i\leq q$, color all $(x-1)/2$ parallel edges in $M_{v_{ij}}$ with
color $j$. Moreover, for vertex $v_x$, color the edge
$v_{ij}v_{(q+1-i)j}$ in $M_{v_x}$ with $j$. The result is a coloring
of $K_x$ with the property that the set of all neighbors
of  every vertex contains $q$ edges of any color, as desired.$\hfill \dashv $ \\
\begin{lemma}\label{lemmaB}
Suppose that $x=tq+r$ is odd and $2\leq r\leq t-1$. There exists
an edge coloring of $K_x$ with $t$ colors such that the set of all
neighbors of every vertex contains at least $q$ edges of any
color.
\end{lemma}
{\it Proof.} Partition the vertices of $K_x$ as $v_1, v_2, \ldots,
v_r$ plus $q$ classes $T_1, \ldots, T_q$ where $T_i$, $1\leq i\leq
q$,  contains $t$ vertices say $v_{i1}, \ldots, v_{it}$ (see $(b)$
of figure \ref{G1}). For each vertex $v_{ij}$  color all $(x-1)/2$
parallel edges in $M_{v_{ij}}$ with color $j$. Moreover, for vertex
$v_r$ (also $v_1$) color the parallel edges in $M_{v_r}$ (also in
$M_{v_1}$) with $1, 2, \ldots, t$ alternatively (also $t, t-1,
\ldots, 1$ alternatively).  Color the remaining edges, i.e. parallel
edges corresponding to $v_2, \ldots, v_{r-1}$ arbitrarily. The
result is a coloring of edges of $K_x$ with the property that for
any
vertex, each color appears on at least $q$ edges, as desired.$\hfill \dashv $ \\

\section{\normalsize{The value of $R_{t-1, t}(K_{1, n})$}}\label{secdelta}

In this section, using  Corollary \ref{cor1}, we determine the
exact value of  $R_{t-1, t}(K_{1, n})$.

\begin{theorem}\label{th3}
Suppose that $x=\left[\frac{nt-1}{t-1}\right]$ and
$q=[\frac{x}{t}]$. Then
$$R_{t-1,
t}(K_{1, n})=\left \{ \begin{array}{cc} x &~~~ {\rm if} ~x=tq+1 ~{\rm{for}} ~x, q ~{\rm{odd}},  \\
 x+1 & {\rm otherwise.}  \end{array}\right. $$
\end{theorem}
{\it Proof.} 
First note that since $x=[\frac{nt-1}{t-1}]$, then
$\frac{nt-1}{t-1}-1<x \leq \frac{nt-1}{t-1}$, or equivalently
$x-x/t+1/t\leq n < x-x/t+1$ and so $n=x-[x/t]=x-q$. If $x$ is
even, then by Corollary \ref{cor1}, $R_{t-1, t}(K_{1, n})=x+1$. So
we may assume that $x$ is odd. We consider three cases as follows.

\vspace{.5cm} \noindent{\bf  Case $1$.} $x=tq+1$, where $q$ is odd.

Consider an edge coloring of $K_x$ with $t$ colors. Suppose first
that any color appears on $q$ edges adjacent to every vertex.
Consider a color $c$, then the subgraph induced by the edges with
color $c$ is $q$-regular and so the sum of degrees of its vertices
is equal to the odd number $xq$, a contradiction. Thus there exist a
vertex $v$ and a color $c$ with the property that $c$ appears on at
most $q-1$ edges adjacent to $v$. Then there are at least
$x-1-(q-1)=x-q=n$ edges adjacent to $v$ such that $c$ does not
appear on these edges. Hence there exists a subgraph $K_{1,n}$
without color $c$ in $K_x$, i.e. $R_{t-1, t}(K_{1, n})\leq x$ and so
by Corollary \ref{cor1}, $R_{t-1, t}(K_{1, n})= x$.

\vspace{.5cm} \noindent{\bf  Case $2$.} $x=tq+1$, where $q$ is even.

In the coloring of $K_x$ given by Lemma \ref{lemmaA}, every
$K_{1,n}$ contains all $t$ colors, i.e. $R_{t-1, t}(K_{1, n})>x$
and so  by Corollary \ref{cor1}, $R_{t-1, t}(K_{1, n})= x+1$.

\vspace{.5cm} \noindent{\bf  Case $3$.} $x=tq+r$, where $2\leq r\leq
t-1$.

In the coloring of $K_x$ given by Lemma \ref{lemmaB}, every
$K_{1,n}$ contains all $t$ colors, i.e. $R_{t-1, t}(K_{1, n})>x$
and so by Corollary \ref{cor1}, $R_{t-1, t}(K_{1, n})= x+1$.$\hfill \dashv $ \\

As a corollary, we have the value of standard Ramsey number
$r_2(K_{1, n})$ (see \cite{sur}).

\begin{corollary}
$r_2(K_{1, n})=2n-\epsilon$  where $\epsilon=1$  if  n  is even
and $\epsilon=0$, otherwise.
\end{corollary}

\section{\normalsize{The value of $R_{t-2, t}(K_{1, n})$}}\label{l-2}

In this section, we determine $R_{t-2, t}(K_{1, n})$. Corollary
\ref{thm3} gives a lower bound and an upper bound for $R_{t-2,
t}(K_{1, n})$ for $t\geq 4$. Let us first settle the case $t=3$.
It is also a special case of multi-color Ramsey numbers for stars
obtained in \cite{ramseyforstars}.

\begin{lemma}\label{lem1,3}
There exists an edges coloring of $K_{3n-2}$ with $3$ colors such
that every vertex contains exactly $n-1$ edges from each color.
\end{lemma}
{\it Proof.} If $3n-2$ is even, then Vizing's Theorem gives a
proper edge coloring of $K_{3n-2}$ with $3n-3$ colors. Divide
these $3n-3$ colors into $3$ new color classes each of which
contains $n-1$ colors to get the desired coloring of $K_{3n-2}$
with $3$ colors. Thus we may assume that $3n-2$ is odd. Then
$K_{3n-2}$ has $3n-2$ matchings each of which contains $(3n-3)/2$
parallel edges. For every vertex, color the corresponding parallel
edges with $1$, $2$ and $3$ respectively to get the desired coloring.$\hfill \dashv $ \\

\begin{theorem}
It holds $R_{1,3}(K_{1, n})=3n-1$.
\end{theorem}
{\it Proof.} Consider an arbitrary edge coloring of $K_{3n-1}$ with
$3$ colors $1$, $2$, $3$ and a vertex $v$. Suppose that $3$ is a
color with maximum number of edges adjacent to $v$. So two colors
$1$ and $2$ appear on at most $2[\frac{3n-2}{3}]$ edges adjacent to
$v$. It is easily seen that $3n-2-2[\frac{3n-2}{3}]\geq n$ and so we
have a $K_{1, n}$ with color $3$, i.e. $R_{1,3}(K_{1, n})\leq 3n-1$.
To prove $R_{1,3}(K_{1, n})\geq 3n-1$, apply Lemma \ref{lem1,3}. In
this coloring of $K_{3n-2}$ every $K_{1, n}$ contains at least $2$
colors and so $R_{1,3}(K_{1, n})>3n-2$.$\hfill \dashv $


For general case $t\geq 4$, we let $R$ stands for $R_{t-2, t}(K_{1,
n})$, $t'=[t/2]$ and $x=[\frac{nt'-1}{t'-1}]$.

\begin{lemma}\label{cliams}
Suppose that $x-2=tq+r$ where $0\leq r\leq t-1$ and $l$ is a natural
number. Then $x-l-2q<n$ iff $t>(2r+4)/l$ when $t$ is even and
$t>1+(2q+2r+4)/l$, otherwise.
\end{lemma}
{\it Proof.} Since $n=x-[\frac{x}{t'}]$, we have $x-l-2q<n$ iff
$[\frac{x}{t'}]<2q+l$ or equivalently, $\frac{x}{t'}<2q+l$. So
$x-l-2q< n$ iff $t>(2r+4)/l$ when $t$ is even and $t>1+(2q+2r+4)/l$,
otherwise. $\hfill \dashv $\\

Theorem \ref{l-2th1}, states the necessary and sufficient conditions
for $R$ being $x+1$.
\begin{theorem}\label{l-2th1}
Suppose that $x-2=tq+r$ where $0\leq r\leq t-1$. Then $R=x+1$ iff
one the following conditions holds.
\begin{itemize}
\item[$(a)$] $r=t-1>2q+4$ and $x$ is even. \item[$(b)$]
$r=t-1>2q+4$ and $x$ and $t$ are odd. \item[$(c)$] $r=t-1$, $x$ is
odd and $t$ and $q+1$ are even. \item [$(d)$] $r<t-2$ and $t>2r+4$
is even. \item [$(e)$] $r<t-2$ and $t>2q+2r+5$ is odd.
\end{itemize}

\end{theorem}
{\it Proof.} We first suppose that $x$ is even and consider three
cases as follows.

\vspace{.5cm} \noindent{\bf  Case $1.1$.} $r=t-1$.

Note that since $x-1=t(q+1)$ is odd, $t$ can't be even. Let
$t>2q+5$. To prove $R=x+1$, using Corollary \ref{thm3}, it is enough
to give a coloring of $K_x$ with $t$ colors such that every $K_{1,
n}$ contains at least $t-1$ colors. By Vizing's Theorem, there
exists a proper edge coloring of $K_x$ with $x-1$ colors. We
partition these $x-1$ colors into $t$ color classes each of which
contains $q+1$ colors to get a coloring of $K_x$ with $t$ colors.
Then every $K_{1, n}$ contains at least $t-1$ colors iff
$x-1-2(q+1)<n$ which holds by the assertion and Lemma \ref{cliams}
for $l=3$. Now let $t\leq 2q+5$. Suppose that an arbitrary edge
coloring of $K_x$ with $t$ colors is given. For each  vertex $v$,
there are least two colors that appear on at most $2(q+1)$ edges of
$E_G(v)$, since $x-1=t(q+1)$. Using Lemma \ref{cliams} for $l=3$, at
least $n$ edges of $E_G(v)$ are colored with the remaining $t-2$
colors, that is, $R\leq x$.

\vspace{.5cm} \noindent{\bf  Case $1.2$.} $r=t-2$.

We now prove $R\neq x+1$ by showing that $R\leq x$. Suppose that an
arbitrary edge coloring of $K_x$ with $t$ colors is given. For each
vertex $v$, there are two colors that appear on at most $2q+1$ edges
of $E_G(v)$, since $x-1=t(q+1)-1$. Using Lemma \ref{cliams} for
$l=2$, at least $n$ edges of $E_G(v)$ are colored with the remaining
$t-2$ colors, that is, there exists a $K_{1, n}$ with at most $t-2$
colors.

\vspace{.5cm} \noindent{\bf  Case $1.3$.} $r<t-2$.

Let either $t>2r+4$ be even or $t>2q+2r+5$ be odd. To prove $R=x+1$,
it is enough to give a coloring of $K_x$ with $t$ colors such that
every $K_{1, n}$ contains at least $t-1$ colors. By Vizing's
Theorem, there is a proper edge coloring of $K_x$ with $x-1$ colors.
We partition these $x-1$ colors into $t-r-1$ color classes each of
which contains $q$ colors plus $r+1$ color classes each of which
contains $(q+1)$ colors to get a coloring of $K_x$ with $t$ colors.
Then every $K_{1, n}$ contains at least $t-1$ colors iff $x-1-2q<n$
which holds by the assertion and Lemma \ref{cliams} for $l=1$, that
is, $R>x$.

Now  suppose that either $t\leq 2r+4$ or $t\leq 2q+2r+5$ is odd.
Suppose that  an arbitrary edge coloring of $K_x$ is given. For each
vertex $v$, there are two colors that appear on at most $2q$ edges
of $E_G(v)$, since $x-1=tq+r+1<t(q+1)-1$. Hence by the assertion and
Lemma \ref{cliams} for $l=1$, at least $n$ edges of $E_G(v)$ are
colored with the remaining $t-2$ colors,  that is, $R\leq x$.

Now suppose that $x$ is odd. We consider three cases as follows.

\vspace{.5cm} \noindent{\bf  Case $2.1$.} $r=t-1$.

Let either $t>2q+5$ be odd or  both of $t$ and $q+1$ be even. We
show that $R=x+1$. Note that since $x-1=t(q+1)$ is even, if $t$ is
odd, then $q+1$ is even. By Lemma \ref{lemmaA}, there exists an edge
coloring of $K_x$ with $t$ colors such that for each vertex $v$,
$E_G(v)$ contains $q+1$ edges of any color. What is left is similar
to the Case $1.1$. If $t\leq 2q+5$ is odd and $q+1$ is even, similar
argument as in the Case $1.1$ yields $R\leq x$. Assume that $q+1$ is
odd and hence $t$ is even. Suppose that an arbitrary edge coloring
of $K_x$ is given. If for each  vertex $v$, $E_G(v)$ contains $q+1$
edges of any color, the induced subgraph on the edges with a fixed
color is $(q+1)$-regular with $x$ vertices, a contradiction. So
there exists a vertex $v$ and a color $c$ such that $E_G(v)$
contains at most $q$ edges with color $c$. So there are two colors
that appear on at most $2q+1$  edges of $E_G(v)$. Since
$x-1-(2q+1)\geq n$, at least $n$ edges of $E_G(v)$ are colored with
the remaining $t-2$ colors, that is, $R\leq x$.

\vspace{.5cm} \noindent{\bf  Case $2.2$.} $r=t-2$.

By the same argument as the Case $1.2$, we get $R\leq x$.

\vspace{.5cm} \noindent{\bf  Case $2.3$.} $r<t-2$.

Let either $t>2r+4$ be even or $t>2q+2r+5$ be odd. By Lemma
\ref{lemmaB}, there exists an edge coloring of $K_x$ such that for
each vertex $v$, $E_G(v)$ contains at least $q$ edges of any color.
What is left is similar to the Case $1.3$.$\hfill \dashv $ \\

Theorem \ref{l-2th2}, states the necessary and sufficient conditions
for $R$ being $x$.

\begin{theorem}\label{l-2th2}
Suppose that $x-2=tq+r$ where $0\leq r\leq t-1$. Then  $R=x$ iff
one the following conditions holds.
\begin{itemize}
\item[$(a)$] $r=t-1$ and $x$ and $q+1$ are odd. \item
[$(b)$] $r<t-2$ and $t\leq 2r+4$ is even. \item [$(c)$] $r<t-2$
and $q+r+3<t\leq 2q+2r+5$ is odd.
\end{itemize}

\end{theorem}
{\it Proof.} Let $p=x-1$, then $p-1=tq+r$. We first suppose that
$p$ is even and consider three cases as follows.

\vspace{.5cm} \noindent{\bf  Case $1.1$.} $r=t-1$.

Let $t$ be even and $q+1$ be odd. By Theorem \ref{l-2th1}, $R\leq
x$. By Vizing's Theorem there exists a proper edge coloring of $K_p$
with $p-1$ colors. We partition these $p-1$ colors into $t-1$
classes each of which contains $q+1$ colors plus a class which
contains $q$ colors to get a coloring of $K_p$ with $t$ colors. Then
every $K_{1, n}$ contains at least $t-1$ colors iff $p-1-(2q+1)<n$
which holds by the assertion and Lemma \ref{cliams} for $l=3$, that
is, $R>p=x-1$ and so $R=x$.

If both of $t$ and $q+1$ are  even then $R>x$ by Theorem
\ref{l-2th1}. Note that the case when  both of $t$ and $q+1$ are odd
is impossible, since $p=t(q+1)$ is even. Assume that $t$ is odd and
$q+1$ is even. If $t>2q+5$, then $R\neq x$ by Theorem \ref{l-2th1}.
Let $t\leq 2q+5$ be odd and $q+1$ be even. Suppose that an arbitrary
edge coloring of $G=K_p$ with $t$ colors is given. For each vertex
$v$, there are  two colors that appear on at most $2q+1$ edges of
$E_G(v)$, since $p-1=t(q+1)-1$. Hence by Lemma \ref{cliams} for
$l=3$, at least $n$ edges  of $E_G(v)$ are colored with the
remaining $t-2$ colors,  that is, $R\leq p=x-1$.

\vspace{.5cm} \noindent{\bf  Case $1.2$.} $r=t-2$.

Suppose that an arbitrary edge coloring of $G=K_p$ with $t$ colors
is given. For each vertex $v$, there are two colors that appear on
at most $2q$ edges of $E_G(v)$, since $p-1=t(q+1)-2$. Hence by Lemma
\ref{cliams} for $l=2$, at least $n$ edges of $E_G(v)$ are colored
with the remaining $t-2$ colors, that is, $R\leq p=x-1$.

\vspace{.5cm} \noindent{\bf  Case $1.3$.} $r<t-2$.

Let either $t\leq 2r+4$ be even or $q+r+3<t\leq 2q+2r+5$ be odd. By
Theorem \ref{l-2th1}, $R\leq x$. By Vizing's Theorem, there exists a
proper edge coloring of $K_p$ with $p-1$ colors. We partition these
$p-1$ colors into $t-r$ color classes each of which contains $q$
colors plus $r$ color classes each of which contains $q+1$ colors to
get a coloring of $K_p$ with $t$ colors. Then every $K_{1, n}$
contains at least $t-1$ colors iff $p-1-2q<n$ which holds by the
assertion and Lemma \ref{cliams} for $l=2$, that is, $R>p=x-1$ and
so $R=x$.

If either $t>2r+4$ is even or $t>2q+2r+5$ is odd, then $R\neq x$ by
Theorem \ref{l-2th1}. Assume that $t\leq q+r+3$ is odd. Suppose that
an arbitrary edge coloring of $G=K_p$ with $t$ colors is given. For
each vertex $v$, there are  two colors that appear on at most $2q$
edges  of $E_G(v)$, since $p-1<t(q+1)-2$. Hence by the assertion and
Lemma \ref{cliams} for $l=2$, at least $n$ edges of $E_G(v)$ are
colored with the remaining $t-2$ colors, that is, $R\leq p=x-1$.

Now suppose that $p$ is odd. We consider three cases as follows.

\vspace{.5cm} \noindent{\bf  Case $2.1$.} $r=t-1$.

So $t$ and $q+1$ are odd. If $t>2q+5$, then by Theorem \ref{l-2th1},
$R=x+1$. Now let $t\leq 2q+5$ be odd. Suppose that an arbitrary edge
coloring of $G=K_p$ with $t$ colors is given. For each vertex $v$,
there are  two colors that appear on at most $2q+1$ edges of
$E_G(v)$, since $p-1=t(q+1)-1$. Hence by Lemma \ref{cliams} for
$l=3$, at least $n$ edges  of $E_G(v)$ are colored with the
remaining $t-2$ colors,  that is, $R\leq p=x-1$.

\vspace{.5cm} \noindent{\bf  Case $2.2$.} $r=t-2$.

By the same arguments as the Case $1.2$, we get $R\leq p=x-1$.

\vspace{.5cm} \noindent{\bf  Case $2.3$.} $r<t-2$.

Let either $t\leq 2r+4$ be even or $q+r+3<t\leq 2q+2r+5$ be odd. If
$r=0$ and $t$ is even, then $t=4$ and so
$x=[\frac{nt'-1}{t'-1}]=2n-1$, which is impossible. By Lemmas
\ref{lemmaA} and \ref{lemmaB}, there exists an edge coloring of
$G=K_p$ with $t$ colors such that for each vertex $v$, $E_G(v)$
contains at least $q$ edges of any color.
What is left is similar to Case $1.3$.$\hfill \dashv$ \\

Theorem \ref{l-2th3}, states the necessary and sufficient conditions
for $R$ being $x-1$.

\begin{theorem}\label{l-2th3}
Suppose that $x-2=tq+r$ where $0\leq r\leq t-1$. Then  $R=x-1$ iff
one the following conditions holds.
\begin{itemize}
\item[$(a)$] $r=1$, $\frac{2q+9}{3}<t\leq q+4$ is odd and $x$
is even \item[$(b)$] $r=1$, $\frac{2q+9}{3}<t\leq q+4$ is odd and
$x$ is odd. \item [$(c)$] $1<r<t-2$ and $\frac{2q+2r+7}{3}<t\leq
q+r+3$ is odd. \item [$(d)$] $r=t-2$ and either $t$ is even or $t>
\frac{2q+2r+7}{3}$ is odd.
\end{itemize}
\end{theorem}

{\it Proof.} Let $p=x-2$, then $p=tq+r$. We first suppose that $x$
is even and consider five cases as follows.

\vspace{.5cm} \noindent{\bf  Case $1.1$.} $r=0$.

If either $t$ is even or $t>2q+5$ is odd, then $R\neq x-1$ by
Theorems \ref{l-2th1} and \ref{l-2th2}. If $q+3<t\leq 2q+5$ is odd,
then $R\neq x-1$, by Theorem \ref{l-2th2}. Now let $t\leq q+3$ be
odd. Note that $q$ is even in this case. Suppose that an arbitrary
edge coloring of $G=K_p$ with $t$ colors is given. For each vertex
$v$, there are  two colors that appear on at most $2q-1$ edges  of
$E_G(v)$, since $p-1=x-3=tq-1$. Hence by the assertion and Lemma
\ref{cliams} for $l=2$, at least $n$ edges of $E_G(v)$ are colored
with the remaining $t-2$ colors, that is, $R\leq p=x-2$ and so by
Corollary \ref{thm3}, $R=x-2$.

\vspace{.5cm} \noindent{\bf  Case $1.2$.} $r=1$.

Let $\frac{2q+9}{3}<t\leq q+4$ be odd. Since $t\leq q+4$, by
Theorems \ref{l-2th1} and \ref{l-2th2}, $R\leq x-1$. By Vizing's
Theorem, there exists a proper edge coloring of $K_p$ with $p-1$
colors. We partition these $p-1$ colors into $t$ color classes each
of which contains $q$ colors to get a coloring of $K_p$ with $t$
colors. Then every $K_{1, n}$ contains at least $t-1$ colors iff
$x-3-2q=p-1-2q<n$ which holds by the assertion and Lemma
\ref{cliams} for $l=3$, that is, $R>p=x-2$ and hence $R=x-1$. If
either $t>q+4$ is odd or $t$ is even, then $R\neq x-1$, by Theorems
\ref{l-2th1} and \ref{l-2th2}. Now let $t\leq \frac{2q+9}{3}$ be
odd. Suppose that an arbitrary edge coloring of $G=K_p$ with $t$
colors is given. For each vertex $v$, there are  two colors that
appear on at most $2q$ edges of $E_G(v)$, since $p-1=x-3=tq$. Hence
by the assertion and Lemma \ref{cliams} for $l=3$, at least $n$
edges of $E_G(v)$ are colored with the remaining $t-2$ colors, that
is, $R\leq p=x-2$.

\vspace{.5cm} \noindent{\bf  Case $1.3$.} $1<r<t-2$.

Let $\frac{2q+2r+7}{3}<t\leq q+r+3$ be odd. Since $t\leq q+r+3$, by
Theorems \ref{l-2th1} and \ref{l-2th2}, $R\leq x-1$. By Vizing's
Theorem, there exists a proper edge coloring of $K_p$ with $p-1$
colors. We partition these $p-1$ colors into $t-r+1$ color classes
each of which contains $q$ colors plus $r-1$ color classes each of
which contains $q+1$ colors to get a coloring of $K_p$ with $t$
colors. Then every $K_{1, n}$ contains at least $t-1$ colors iff
$x-3-2q=p-1-2q<n$ which holds by the assertion and Lemma
\ref{cliams} for $l=3$, that is, $R>p=x-2$ and so $R=x-1$.

If either $t>q+r+3$ is odd or $t$ is even, then $R\neq x-1$ by
Theorems \ref{l-2th1} and \ref{l-2th2}. Now let $t\leq
\frac{2q+2r+7}{3}$ be  odd. Suppose that an arbitrary edge coloring
of $G=K_p$ with $t$ colors is given. For each vertex $v$, there are
two colors that appear on at most $2q$ edges  of $E_G(v)$, since
$p-1=x-3=tq+r-1<t(q+1)-3$. Hence by the assertion and Lemma
\ref{cliams} for $l=3$, at least $n$ edges  of $E_G(v)$ are colored
with the remaining $t-2$ colors, that is, $R\leq p=x-2$.

\vspace{.5cm} \noindent{\bf  Case $1.4$.} $r=t-2$.

Let either $t$ be even or $t>\frac{2q+2r+7}{3}$ be odd.  By Theorems
\ref{l-2th1} and \ref{l-2th2}, $R\leq x-1$.  By Vizing's Theorem,
there exists a proper edge coloring of $K_p$ with $p-1$ colors. We
partition these $p-1$ colors into $t-3$ color classes each of which
contains $q+1$ colors plus $3$ color classes each of which contains
$q$ colors to get a coloring of $K_p$ with $t$ colors. Then every
$K_{1, n}$ contains at least $t-1$ colors iff $x-3-2q=p-1-2q<n$
which holds by the assertion and Lemma \ref{cliams} for $l=3$, that
is, $R>p=x-2$. Therefore $R=x-1$.

 Now let $t\leq \frac{2q+2r+7}{3}$ be
odd. Suppose that an arbitrary edge coloring of $G=K_p$ with $t$
colors is given. For each vertex $v$, there are two colors that
appear on at most $2q$ edges  of $E_G(v)$, since $p-1=x-3=t(q+1)-3$.
Hence by the assertion and Lemma \ref{cliams} for $l=3$, at least
$n$ edges  of $E_G(v)$ are colored with the remaining $t-2$ colors,
that is, $R\leq p=x-2$.

\vspace{.5cm} \noindent{\bf  Case $1.5$.} $r=t-1$.

Hence $t$ is odd. If $t>2q+5$, then by Theorem \ref{l-2th1}, $R\neq
x-1$. Now let $t\leq 2q+5$ be odd. Using Lemma \ref{cliams} for
$l=3$, we have $x-3-2q\geq n$ and so for each edge coloring of
$G=K_p$ with $t$ colors, $n$ edges of $E_G(v)$ are colored with at
most $t-2$ colors, that is $R\leq p=x-2$.

 Now suppose that $x$ is odd. We
consider five cases as follows.

\vspace{.5cm} \noindent{\bf  Case $2.1$.} $r=0$.

The proof is similar to the Case $1.1$. Note that when $t$ is odd,
$q$ can't be even.

\vspace{.5cm} \noindent{\bf  Case $2.2$.} $r=1$.

Let $\frac{2q+9}{3}<t\leq q+4$ be odd and $q$ be even. By Lemma
\ref{lemmaA}, there exists an edge coloring of $K_p$ with $t$ colors
such that every vertex contains  $q$ edges of any color. What is
left is similar to Case $1.2$. Note that the case when both of $t$
and $q+1$ are odd is impossible.

\vspace{.5cm} \noindent{\bf  Case $2.3$.} $1<r<t-2$.

Let $\frac{2q+2r+7}{3}<t\leq q+r+3$ be odd. By Lemma \ref{lemmaB},
there exists an edge coloring of $K_p$ with $t$ colors such that
every vertex contains at least $q$ edges of any color. What is left
is similar to  Case $1.3$.

\vspace{.5cm} \noindent{\bf  Case $2.4$.} $r=t-2$.

Let either $t$ be even or $t>\frac{2q+2r+7}{3}$ be odd. By Lemma
\ref{lemmaB}, there exists an edge coloring of $K_{p}$ with $t$
colors such that every vertex contains at least $q$ edges of any
color. What is left is similar to  Case $1.4$.

\vspace{.5cm} \noindent{\bf  Case $2.5$.} $r=t-1$.

If either $t$ is even or $t>2q+5$ is odd, then $R\neq x-1$ by
Theorems \ref{l-2th1} and \ref{l-2th2}. What is left is similar to
the Case $1.5$.$\hfill \dashv $

\begin{corollary}
$R=x-2$ iff none of the conditions stated in Theorems
\ref{l-2th1}, \ref{l-2th2} and \ref{l-2th3} holds.
\end{corollary}

\vspace{1.5cm}


\footnotesize

\begin{thebibliography}{99}


\bibitem{ramseyforstars}
{S.A. Burr and J.A. Roberts}, On Ramsey numbers for stars,
{\it Utilitas. Math.} 4 (1973) 217-220.

\bibitem{chli1}
{ K.M. Chung, M.L. Chung and C.L. Liu}, A generalization of
Ramsey theory for graphs-with stars and complete graphs as
forbidden subgraphs, {\it Congr. Numer.} 19 (1977) 155-161.

\bibitem{chli2}
{ K.M. Chung and C.L. Liu}, A generalization of Ramsey theory
for graphs, {\it Discrete Math.} 2 (1978) 117-127.



\bibitem{gerencser}
{ L. Gerencs\'{e}r and A. Gy\'{a}rf\'{a}s}, On Ramsey-type
problems, {\it Ann. Univ. Sci. Budapest E\"{o}tv\"{o}s.} 10 (1967)
167-170.


\bibitem{gyarfas}
{ A. Gy\'{a}rf\'{a}s, G.N. S\'{a}rk\"{o}zy and S. Selkow},
Coverings by few monochromatic pieces - a
transition between two Ramsey problems,
manuscript submitted in 2011.


\bibitem{harborth}
{ H. Harborth and M. M$\ddot{o}$ller}, Weakened Ramsey numbers,
{\it Discrete Applied Math.} 95 (1999) 279-284.

\bibitem{khomidiforests}
{ A. Khamseh and R. Omidi}, A generalization of Ramsey theory for
linear forests, {\it Int. J. Comput. Math.} 89 (2012) 1303-1310.

\bibitem{meenakshi}
{ R. Meenakshi and  P.S. Sundararaghavan}, Generalized Ramsey
numbers for paths in $2$-chromatic graphs, {\it Internat. J. Math.
Sci.} 9 (1986) 273-276.

\bibitem{sur}
{ S.P. Radziszowski}, Small Ramsey numbers, {\it Electronic J.
Combin.} 1 (1994) Dynamic Surveys, DS1.13 (August 22, 2011).













\end{thebibliography}
 \end{document}